\documentclass[twoside,10pt]{amsart}
\usepackage{graphicx,amsmath,amsfonts,amsthm,amssymb,fancyhdr,lscape, pdflscape}
\newtheorem{thm}{Theorem}[section]

\newtheorem{lem}[thm]{Lemma}

\theoremstyle{definition}
	
\theoremstyle{remark}

\def\beq{\begin{eqnarray}}
\def\eeq{\end{eqnarray}}
\def\bsp{\begin{split}}	
\def\esp{\end{split}}

\newcommand{\ov}[3]{ {#1}^{(#2)}_{#3} }
\newcommand{\ox}[3]{ {#1}^{[#2]}_{#3} }

\newcommand{\be}{\begin{equation}}
\newcommand{\ee}{\end{equation}}

\begin{document}

\begin{abstract}

We study the existence of a non-spacelike isometry, $\zeta$, in higher dimensional Kundt spacetimes with constant scalar curvature invariants ($CSI$). We present the particular forms for the null or timelike Killing vectors and a set of constraints for the metric functions in each case.  Within the class of $N$ dimensional $CSI$ Kundt spacetimes, admitting a non-spacelike isometry, we determine which of these can admit a covariantly constant null vector that also satisfy $\zeta_{[a;b]} = 0$.

\end{abstract}

\title{Killing vectors in higher dimensional Spacetimes \\ with constant scalar curvature invariants}

\author{\textbf{David McNutt, Nicos Pelavas, Alan Coley}}
\address{Department of Mathematics and Statistics, Dalhousie University, Halifax, Nova Scotia, Canada B3H 3J5}
\email{aac, mcnuttd, pelavas@mathstat.dal.ca}
\date{\today}


\maketitle

\section*{Introduction}

An $N$ dimensional differentiable manifold of Lorentzian signature for which all polynomial scalar curvature invariants constructed from the Riemann tensor and its covariants derivatives are constant is called a $CSI$ spacetime.  There are many examples of $CSI$ spacetimes in general relativity and other gravity theories.  The higher dimensional pp-wave spacetimes, which are exact solutions of supergravity and string theory of Ricci type N, have vanishing polynomial scalar curvature invariants. We call such spacetimes $VSI$ and note that the set of all $VSI$ spacetimes is a subset of the set of all $CSI$ spacetimes. 

There are further subdivisions in the set of all $CSI$ spacetimes depending on distinguishing properties of the spacetimes:

\begin{itemize}  

\item ${CSI_R}$ - The set of  all reducible $CSI$ spacetimes that can be built from $VSI$ and $H$ by (i) warped products (ii) fibered products, and (iii) tensor sums. 
\item ${CSI_F}$ - All spacetimes for which there exists a frame with a null vector $\ell$ such that all components of the Riemann tensor and its covariants derivatives in this frame have the property that (i) all positive boost weight components (with respect to $\ell$) are zero and (ii) all zero boost weight components are constant. 
\item ${CSI_K}$ -  Those $CSI$ spacetimes that belong to the (higher dimensional) Kundt class; the so-called Kundt $CSI$ spacetimes.

\end{itemize}

For a Riemannian manifold every $CSI$ spacetime is homogeneous; this is not true for Lorentzian manifolds.  However, for every $CSI$ spacetime with particular constant invariants there is a homogeneous spacetime (not necessarily unique) with precisely the same constant invariants.  This suggests that $CSI$ spacetimes can be constructed from $H$ and $VSI$ (e.g., ${CSI_R}$). In particular, the relationship between ${CSI_R}$, ${CSI_F}$, ${CSI_K}$ and especially with $CSI \backslash H$ were studied in arbitrary dimensions \cite{CSI} (and considered in more detail in the four dimensional case). We note that by construction ${CSI_R}$ is at least of Weyl type $II$ (i.e., of type $II$, $III$,  $N$ or $O$ \cite{class}), and by definition ${CSI_F}$ and ${CSI_K}$ are at least of Weyl type $II$ (more precisely, at least of Riemann type $II$). In four dimensions ${CSI_R}$, ${CSI_F}$ and ${CSI_K}$ are closely related and that if a spacetime was $CSI$ then it is either homogeneous or belong to the Kundt $CSI$ spacetimes. (\cite{4DCSI} \cite{4DKundt}, \cite{CSI}). It is conjectured that in higher dimensions this holds true as well. 

It was shown in \cite{VSISugra} that more generally all Ricci type N $VSI$ spacetimes and  some Ricci type III $VSI$ spacetimes (assuming appropriate sources) are solutions to type IIB supergravity; it was argued that these are also solutions to other supergravity theories as well. There are many supergravity $CSI$ spacetimes \cite{CSISugra}. The $CSI$ spacetime $AdS_d \times S^{D-d}$ for fixed $D$ is a solution to supergravity as it preserves the maximal number of superymmetries. There are other $CSI$ solutions which admit supersymmetries: the $AdS$ gyratons (\cite{AdSGyr0}, \cite{AdSGyr1}) or chiral null models \cite{ChiralNull}. However, the full set of $CSI$ spacetimes which admit supersymmetries has not been studied in detail. 

In \cite{VSISugra} it was noted that a particular solution of type IIB supergravity  admit a non-spacelike isometry in order to admit a supersymmetry. Furthermore in \cite{VSISugra} it was proven that the only $VSI$ spacetimes which admit a null or timelike Killing vector are those which already admit a covariantly constant null vector $\ell = \frac{\partial}{\partial{v}}$. This implies the set of $VSI$ spacetimes which satisfy type IIB supergravity belong to the subset of $VSI$ Kundt spacetimes \cite{Kundt} with no $v$ dependence. We intend to generalize this result to the set of $CSI$ Kundt spacetimes \cite{Kundt}. That is, we shall determine the set of $CSI$ Kundt spacetimes admitting a null or timelike Killing vector. 

\subsection*{Kundt Spacetimes}

In higher dimensions it was shown that $\ell$ is geodesic, non-expanding, shear-free and non-twisting in $VSI$  spacetimes \cite{Higher}; its covariant derivative takes the form:

\beq \ell_{a ; b} = L_{11} \ell_a \ell_b  + L_{1i} \ell_a m^i_{~b} + L_{i1} m^i_{~a} \ell_b. \label{covariantl} \eeq   

\noindent For locally homogeneous spacetimes, in general there exists a null frame in which the $L_{ij}$ are constants. We therefore anticipate that in $CSI$ spacetimes that are not locally homogeneous $L_{ij}=0$. For higher dimensional $CSI$ spacetimes with $L_{ij}=0$, the Ricci and Bianchi identities appear to be identically satisfied \cite{bianchi}. A higher dimensional spacetime admitting a null vector $\ell$ which is geodesic, non-expanding, shear-free and non-twisting, will be denoted as a higher dimensional Kundt spacetime. 

It was shown in \cite{CSI} that there exists a local coordinate system $(u,v,x^i)$ such that 

\beq  d s^2=2du\left(d v+H(v,u,x^k) du+W_{ i}(v,u,x^k)d x^i\right)+\tilde{g}_{ij}(u,x^k)dx^id x^j \label{Kundt} \eeq

\noindent and the only coordinate transformations preserve the Kundt form \cite{CSI} are the following:

\beq &(v',u',x'^i) =(v,u, f^e(x^g)),~~with~J^e_{~f} = \frac{\partial f^e}{\partial x^f}~~~~~~~~~~~~~&  \nonumber \\ 
&H' = H,~W'_e = W_f (J^{-1})^f_{~e},~ g'_{ef} = g_{gh}(J^{-1})^g_{~e}(J^{-1})^h_{~f},& \label{type1ct} \\
& &\nonumber \\
&(v',u',x'^i) = (v+h(u,x^e), u, x^i)~~~~~~~~~~~~~~~~~~~~~~~~~~~~& \nonumber \\
& H' = H-h_{,u},~~W'_e = W_e - h_{,e},~~g'_{ef} = g_{ef}, & \label{type2ct} \\
& &\nonumber \\
&(v', u', x'^i) = (v/g_{,u}, g(u), x^i)~~~~~~~~~~~~~~~~~~~~~~~~~~~~~~~~~& \nonumber \\
& H' = \frac{1}{ g^2_{,u} }(H + v \frac{ g_{,uu} }{ g_{,u} } ),~~W'_e = \frac{1}{ g_{,u} } W_e,~~ g'_{ef} = g_{ef}.& \label{type3ct} \eeq

\noindent  Furthermore it was shown that for all $CSI \subset {CSI_F}  \bigcap {CSI_K}$, there always exists (locally) a coordinate transformation $(v',u',x'^i)=(v,u,f^i(u;x^k))$ preserving the form of the Kundt metric, and such that 

\beq \tilde{g}_{ij}\equiv \tilde{g}'_{kl}\frac{\partial f^k}{\partial x^i}\frac{\partial f^l}{\partial x^j}, \quad
\tilde{g}'_{ij,u'}=0 \nonumber \eeq 

\noindent where $d S_H^2=\tilde{g}'_{ij}dx'^id x'^j$ is a locally homogeneous space.

\subsection*{ {$CSI$}$_0$ Kundt spacetimes}

For a particular coframe, 

\beq &n = dv + H(u,v,x^3) du + w_e(u,v,x^e) dx^e,~~\ell = du,~~m^i = m^i_{~e} dx^e& \nonumber \eeq 

\noindent with $m^i_{e}m_{if} = g_{ef}$, the non-zero frame connection components are:

\beq & \Gamma_{21i} = \frac{D_1 W_i}{2},~~\Gamma_{212} = D_1 H,~~ \Gamma_{2i2} = D_i H - D_2 W_i, &  \label{csigamma1} \\
& \Gamma_{i12} = \frac{D_1 W_i}{2},~~\Gamma_{i21} = \frac{D_1W_i}{2},~~\Gamma_{i2j} = \frac{A_{ij}}{2},~~\Gamma_{ij2} = \frac{A_{ij}}{2}, &  \label{csigamma2} \\ 
&\Gamma_{ijk} = -\frac12 \left( D_{ijk} + D_{jki} - D_{kij}  \right) &.  \label{csigamma3} \eeq 

\noindent Where the tensors involved are  written in terms of $m_{ie}$ and it's inverse $D_{ijk} = 2 m_{ie,f}m_{[j}^{~~e}m_{k]}^{~~f}$, and $A_{ij} = D_{[j}W_{i]} - D_{kji}W^k = 2W_{[i;j]}$.

The linearly independent components of the Riemann tensor with boost weight 1 and 0 may be written as: 

\beq R_{ 1 2 1 i}&=&-\frac 12 W_{ i,vv} \nonumber \\
R_{ 1 2 1 2}&=& -H_{,vv}+\frac 14\left(W_{{i},v}\right)\left(W^{ i,v}\right), \nonumber \\
R_{ 1 2 i j}&=& W_{[ i}W_{ j],vv}+W_{[ i; j],v}, \nonumber \\
R_{ 1 i 2 j}&=& \frac 12\left[-W_{ j}W_{ i,vv}+W_{ i; j,v}-\frac 12 \left(W_{  i,v}\right) \left(W_{ j,v}\right)\right], \nonumber \\
R_{ i j \hat i \hat j}&=&\tilde{R}_{ i j \hat i \hat j}. \nonumber \eeq 

\noindent The spacetime will be $CSI_0$ if there exists a frame $\{ \ell, n, m^i \}$, a constant $\sigma$, anti-symmetric matrix ${\sf a}_{\hat i\hat j}$, and symmetric matrix ${\sf s}_{\hat i\hat j}$ such that: 

\beq W_{\hat i,vv} &=& 0, \label{Wcsi1}\\
H_{,vv}-\frac 14\left(W_{\hat i,v}\right)\left(W^{\hat i,v}\right) &=& \sigma, \label{HWiSigma}\\
W_{[\hat i;\hat j],v} &=& {\sf a}_{\hat i\hat j}, \\
W_{(\hat i;\hat j),v}-\frac 12 \left(W_{\hat i,v}\right)\left(W_{\hat j,v}\right) &=& {\sf s}_{\hat i\hat j}, \label{Wcsi4}\eeq 

\noindent and the components $\tilde{R}_{ i j \hat i \hat j}$ are all constants (i.e., $d S^2_H$ is curvature homogeneous). We note that \eqref{Wcsi1} and \eqref{HWiSigma} imply that the metric functions take the following form: 

\beq W_{ i}(v,u,x^e)=v{W}_{ i}^{(1)}(u,x^e)+{W}_{ i}^{(0)}(u,x^e), \label{Wiform} \\ H(v,u,x^e)=\frac{v^2}{8}\left[4\sigma+({W}_i^{(1)})({W}^{(1)i})\right]+v{H}^{(1)}(u,x^e)+{H}^{(0)}(u,x^e). \label{Hform} \eeq

\section*{The Killing Equations} 

Let $\zeta = \zeta_1 n + \zeta_2 \ell + \zeta_i m^i$ be a Killing vector field in a $CSI$ Kundt spacetime; it satisfies the Killing equations for $a,b \in [1,N]$ 

\beq \zeta_{a,b} + \zeta_{b,a} - 2\Gamma^c_{~(ab)}\zeta_c= 0. \nonumber \eeq

\noindent To simplify the analysis of these equations, we choose new coordinates where one of the Killing vectors of the transverse space, $Y$, has been rectified so that locally it behaves as a translation; i.e., $Y = A \frac{\partial}{\partial x^3}$. In this coordinate system $g_{33}$ will be constant, and so it is possible to pick a coframe with an upper-triangular matrix $m^i_{~e}$ and $m^3_{~3}$ constant \cite{Gugg}. This choice of coframe causes $\Gamma_{3ij}$ and $\Gamma_{3(ij)}$ $ \forall i,j \in [3,N]$ to vanish.  Rotating the frame so that the spatial component of $\zeta$ is locally aligned with $m^3$, $\zeta$ takes the form $ \zeta = \zeta_1 n + \zeta_2 \ell + \zeta_3 m^3 $

The components $\zeta_1$ and $\zeta_3$ may be partially integrated from the equations with indices $(11)$, $(13)$, $(3i)$:

\beq & \zeta_1 = ~~\zeta_1(u,x^3),~~ \zeta_3 = -D_3( \zeta_1)v + \ov{\zeta}{0}{3}(u,x^e), & \label{kvcomp13} \eeq

\noindent where $\ov{\zeta}{0}{3}$ satisfies the following differential equations from $(3i)$:

\beq D_i \ov{\zeta}{0}{3} + \ov{W}{0}{i} D_3( \zeta_1) = 0,~~~ D_iD_3\zeta_1 - \ov{W}{1}{i} D_3 \zeta_1 = 0 \label{zeta3iskey} \eeq

\noindent The tensors $\Gamma_{2i2} = D_i H - D_2 W_i$ and $A_{mn} = D_{[n}W_{m]}$ may be expanded into orders of $v$:

\beq \Gamma_{2i2} &=& \overbrace{(D_i \sigma^* - \sigma^* \ov{W}{1}{i})}^{ \ov{\Gamma}{2}{i} } \frac{v^2}{8} + \overbrace{(D_i \ov{H}{1}{~} - \frac14 \ov{W}{0}{i} \sigma^* - D_2\ov{W}{1}{i})}^{ \ov{\Gamma}{1}{i} } v \nonumber \\
 &~& + \overbrace{ D_i \ov{H}{0}{~}- \ov{W}{0}{i} \ov{H}{1}{~} - D_2 \ov{W}{0}{i} + \ov{H}{0}{~} \ov{W}{1}{i} }^{ \ov{\Gamma}{0}{i} },  \label{note1} \\
A_{ij} &=& \overbrace{ 2D_{[j} \ov{W}{1}{i]} }^{ \ov{A}{1}{{ij}} } v + \overbrace{ 2D_{[j} \ov{W}{0}{i]} - 2\ov{W}{0}{[j} \ov{W}{1}{i]} }^{\ov{A}{0}{{ij}} }. \nonumber \eeq

\noindent Here $\sigma^* \equiv 4\sigma + \ov{W}{1}{i} W^{(1)i}$ and the metric functions $H$ and $W_i = m_i^{~e} w_e$ are of the form \eqref{Wiform} and \eqref{Hform}. Substituting these into the equation with indices $(21)$ yields $\zeta_2$ in orders of $v$:

\beq & \zeta_2 = \overbrace{ ( \frac{ \sigma^* \zeta_1}{4} - \ov{W}{1}{3} D_3( \zeta_1) ) }^{\ov{\zeta}{2}{2} } \frac{v^2}{2}   + \overbrace{ ( \ov{W}{1}{3} \ov{\zeta}{0}{3} - D_2 \zeta_1 + \ov{H}{1}{~} \zeta_1 ) }^{\ov{\zeta}{1}{2} } v + \ov{\zeta}{0}{2}(u,x^e). & \label{kvcomp2} \eeq 

Our primary interest are those $CSI$ spacetimes which do not admit covariantly constant null vectors, since the existence of Killing vectors in $CCNV$ spacetimes was considered in \cite{CCNVmcnutt}. The analysis will be restricted to non-spacelike Killing vectors, $|\zeta| \leq 0$. Using the definition of the vector components given above the magnitude is expanded into orders of $v$: 

\beq & \frac{-\sigma^*}{4} (\zeta_1)^2 + \ov{W}{1}{3} D_3( \zeta_1) \zeta_1 + (D_3( \zeta_1))^2 \leq 0 & \label{magv2}\\
&  \zeta_1 (\ov{W}{1}{3} \ov{\zeta}{0}{3} - D_2 \zeta_1 + \ov{H}{1}{~} \zeta_1) + D_3( \zeta_1) \ov{\zeta}{0}{3} = 0 & \label{magv1}\\
& (\ov{\zeta}{0}{3})^2 - 2 \zeta_1 \ov{\zeta}{0}{2} \leq 0. & \label{magv0} \eeq



\noindent The remaining Killing equations, with indices $22$, $23$ and $2n$ are now expanded into orders of $v$, giving the following set of equations: 

\beq &\ov{\Gamma}{2}{3} D_3( \zeta_1) = 0, & \label{3ke22}   \\
& D_2 \ov{\zeta}{2}{2} + \frac14 \sigma^* \ov{\zeta}{1}{2} - \ov{H}{1}{~} \ov{\zeta}{2}{2} - \frac14 \ov{\Gamma}{1}{3} D_3( \zeta_1) + \frac14 \ov{\Gamma}{2}{3} \ov{\zeta}{0}{3} = 0, & \label{2ke22}   \\
& D_2 \ov{\zeta}{1}{2} + \frac14 \sigma^* \ov{\zeta}{0}{2} - \ov{H}{0}{~} \ov{\zeta}{2}{2} - \ov{\Gamma}{0}{3} D_3( \zeta_1) + \ov{\Gamma}{1}{3} \ov{\zeta}{0}{3} = 0,  & \label{1ke22}   \\ 
&D_2 \ov{\zeta}{0}{2} - \ov{H}{0}{~} \ov{\zeta}{1}{2} + \ov{H}{1}{~} \ov{\zeta}{0}{2} + \ov{\Gamma}{0}{3} \ov{\zeta}{0}{3} = 0, &  \label{0ke22} \eeq

\beq & \frac14 \sigma^* D_3( \zeta_1) + D_3 \ov{\zeta}{2}{2} - \ov{W}{1}{3} \ov{\zeta}{2}{2} - \frac14 \ov{\Gamma}{2}{3} \zeta_1 = 0, &   \label{2ke23} \\
& D_2 D_3( \zeta_1) - \ov{H}{1}{~} D_3( \zeta_1) - D_3 \ov{\zeta}{1}{2} + \ov{W}{0}{3} \ov{\zeta}{2}{2} + \ov{\Gamma}{1}{3} \zeta_1 = 0, &   \label{1ke23} \\
&D_2 \ov{\zeta}{0}{3} + \ov{H}{0}{~} D_3( \zeta_1) + D_3 \ov{\zeta}{0}{2} - \ov{W}{0}{3} \ov{\zeta}{1}{2} - \ov{\Gamma}{0}{3} \zeta_1 + \ov{W}{1}{3} \ov{\zeta}{0}{2} = 0, & \label{0ke23}    \eeq

\beq & D_n \ov{\zeta}{2}{2} - \ov{W}{1}{n} \ov{\zeta}{2}{2} - \frac14 \ov{\Gamma}{2}{n} \zeta_1 +  \ov{A}{1}{{3n}} D_3( \zeta_1) = 0, &  \label{2ke2n} \\
& D_n \ov{\zeta}{1}{2} - \ov{W}{0}{n} \ov{\zeta}{2}{2} - \ov{\Gamma}{1}{n} \zeta_1 - \ov{A}{1}{{3n}} \ov{\zeta}{0}{3} + \ov{A}{0}{{3n}} D_3( \zeta_1) = 0, & \label{1ke2n} \\
& D_n \ov{\zeta}{0}{2} - \ov{W}{0}{n} \ov{\zeta}{1}{2} - \ov{\Gamma}{0}{n} \zeta_1 + \ov{W}{1}{n} \ov{\zeta}{0}{2} - \ov{A}{0}{{3n}} \ov{\zeta}{0}{3} = 0. &  \label{0ke2n} \eeq

\noindent The analysis splits into subcases arising from \eqref{3ke22} where either $D_3( \zeta_1)$ or $\ov{\Gamma}{2}{3}$ are assumed to vanish seperately. 

\section*{ Implications of $\zeta_{[a;b]} = 0$}

Before each case is analyzed it will be beneficial to examine the anti symmetrization of $\zeta_{a;b} = 0$ to determine the set of $CSI$ spacetimes admitting a covariantly constant non-spacelike vector. Non-spacelike Killing vectors in $CCNV$ $CSI$ spacetimes has already been studied in \cite{CCNVmcnutt} as such if a $CSI$ spacetime is shown to be $CCNV$ it may be disregarded in the current analysis. Conversely it is of interest to determine when a $CSI$ spacetime  admits a Killing vector but cannot admit a covariantly constant vector. 

Using the form of $\zeta$ given above, the vanishing of $\zeta_{[a;b]} $, yields the following equations:


\beq D_2 \zeta_1 - D_1 \zeta_2 - \Gamma^1_{~12} \zeta_1 = 0, \\
D_3 \zeta_1 - D_1 \zeta_3 - 2 \Gamma^1_{~[13]} \zeta_1 = 0, \label{azeta13} \\
2 \Gamma^1_{~[1n]} \zeta_1 = 0, \label{azeta1n} \\
D_3 \zeta_2 - D_2 \zeta_3 + \Gamma^1_{~32} \zeta_1 = 0,  \\
D_n \zeta_2 +  \Gamma^1_{~n2} \zeta_1 = 0,  \\
D_n \zeta_3 - 2 \Gamma^1_{~[3n]} \zeta_1  = 0,  \\
\Gamma^1_{~[nm]}\zeta_1  = 0 \label{azetanm}. \eeq

\noindent  Assuming $\zeta_1 \neq 0 $ and expanding \eqref{azeta1n} implies $\Gamma^1_{~[1n]} = \ov{W}{1}{n} = 0$. Similarly $2\Gamma^1_{~[13]} = \ov{W}{1}{3}$ and so equation  \eqref{azeta13} gives $\ov{W}{1}{3} = 2 D_3 ln (\zeta_1)$.  Equation \eqref{azetanm} implies that $A_{nm}$ must vanish. Using \eqref{note1} we may summarize these observations as

\begin{lem}\label{lem:azeta0}

For those spacetimes admitting a vector $\zeta$  such that $\zeta_{(a;b)} = 0$ and $\zeta_{[a;b]} = 0$ it is necessary that the metric functions $W_i$ satisfy the following: 
 
\beq & \ov{W}{1}{3} = 2D_3 ln (\zeta_1),~~ \ov{W}{1}{n} = 0,& \nonumber \\
& A_{nm} = 2D_{[m}\ov{W}{0}{n]} = 0. &\nonumber  \eeq 

\end{lem}

\noindent The remaining equations are: 

\beq D_2 \zeta_1 - D_1 \zeta_2 - \Gamma_{212} \zeta_1 = 0, \label{azeta12} \\ 
D_3 \zeta_2 - D_2 \zeta_3 + \Gamma_{232} \zeta_1 = 0, \label{azeta23} \\
D_n \zeta_2 +  \Gamma_{2n2} \zeta_1 = 0, \label{azeta2n} \\
D_n \zeta_3- A_{3n} \zeta_1 = 0. \label{azeta3n} \eeq

\noindent These will be studied once the  analysis of the Killing equations has been completed. 

\section*{Case 1: $D_3( \zeta_1)=0$}

\noindent Setting $D_3( \zeta_1)$ equal to zero we obtain	

\beq &\zeta_1 = ~~\ox{\zeta}{0}{1}(u),~~\zeta_3 = \ov{\zeta}{0}{3}(u)&   \\ 
&\zeta_2 = \overbrace{ ( \frac{ \sigma^* \zeta_1 }{4}  ) }^{\ov{\zeta}{2}{2} } \frac{v^2}{2} + \overbrace{ ( \ov{W}{1}{3} \zeta_3 - D_2 \zeta_1 + \ov{H}{1}{~} \zeta_1 ) }^{\ov{\zeta}{1}{2} } v + \ov{\zeta}{0}{2}(u,x^e).&  \label{c1zeta2} \eeq

\noindent The non-spacelike conditions are now

\beq & -\sigma^* (\zeta_1)^2  \leq 0,~~ \zeta_1 (\ov{W}{1}{3} \zeta_3 - D_2 \zeta_1 + \ov{H}{1}{~} \zeta_1) = 0,~~ (\zeta_3)^2 - \zeta_1 \ov{\zeta}{0}{2} \leq 0 & \eeq

\noindent so either $\zeta_1$ vanishes and $\zeta$ is a null Killing vector or $\zeta_1 \neq 0$ and $\sigma^* \geq 0$. 

\subsection*{Case 1.1 : $\zeta_1 = 0$ }

\noindent If $\zeta_1$ is allowed to vanish, the remaining non-spacelike conditions imply that $\zeta_3 = 0$ and so the Killing vector is of the form $\zeta = \zeta_2 \ell$. In light of the special form of $\zeta_2$ it must be a function of only $u$ and the spatial coordinates $x^e$. The remaining Killing equations are 

\beq \sigma^* \ov{\zeta}{0}{2} = 0 \label{c11sigma} \\
D_2 \ov{\zeta}{0}{2} + \ov{H}{1}{~} \ov{\zeta}{0}{2} = 0 \label{c11H1eqn} \\
D_3 \ov{\zeta}{0}{2} + \ov{W}{1}{3} \ov{\zeta}{0}{2} = 0 \\
D_n \ov{\zeta}{0}{2} + \ov{W}{1}{n} \ov{\zeta}{0}{2} = 0. \eeq

\noindent The vanishing of $\sigma^*$ in the first term \eqref{c11sigma} implies $\ov{W}{1}{i} W^{(1)i} = - 4 \sigma$ where $\ov{W}{1}{i} = m_i^{~e}\ov{w}{1}{e}$ and hence 

\beq \ov{W}{1}{i} W^{(1)i} = g^{ef}\ov{W}{1}{e} \ov{W}{1}{f} = - 4 \sigma. \label{c11metric} \eeq

\noindent Since the transverse metric is Riemannian, it is positive-definite and restricts the value of $\sigma$ to be less than or equal to zero. 

\subsection*{Case 1.1.1:}

If $\sigma = 0$, this implies $\ov{W}{1}{e} = 0$ for all $e \in [3,N]$. The vector component $\zeta_2$ will be a function of $u$ only and the remaining equation \eqref{c11H1eqn} determines the  metric function 

\beq \ov{H}{1}{~}(u) = - D_2 ln(\zeta_2). \label{c111H1} \eeq

\noindent One may always make a coordinate transform of the form \eqref{type3ct} to set $\ov{H}{1}{~} = 0$,  so the metric is independent of the null coordinate $v$  and $\zeta = \ell = \frac{\partial}{\partial{v}}$, implying that $\zeta$ is a covariantly constant null vector.     

\subsection*{Case 1.1.2: }
 
If $\sigma < 0$, one may solve for the metric functions $\ov{W}{1}{i}$ and $\ov{H}{1}{~}$ in terms of $\zeta(u,x^e)$: 

\beq \ov{H}{1}{~}(u,x^e) = - D_2 ln(\zeta_2),~~ \ov{W}{1}{i}(u,x^e) = -D_i ln(\zeta_2). \nonumber \eeq

\noindent These $CSI$ spacetimes do not admit a covariantly constant vector. To see this, assume $\sigma < 0 $ and consider equations  \eqref{azeta12} - \eqref{azeta2n}; the first two are automatically satsified while the last implies that $\zeta_2$ is a function of $u$ only. This forces the $\ov{W}{1}{i}$ to all vanish, leading to the contradiction: $0 = \sigma^* = \sigma <0$, hence these spacetimes do not admit a $CCNV$.

Given a null vector of the form, $\zeta = \zeta(u,x^e) \ell$, it will be a Killing vector for the $CSI$ spacetime with a locally homogeneous transverse space and metric functions:

\beq &H = -(ln \zeta)_{,u} v + \ov{H}{0}{~}(u, x^e),~~ W_e = -(ln \zeta)_{,e} v + \ov{W}{0}{e}(u,x^f) .& \label{c112metricfn} \eeq

\noindent The vanishing of the function $\sigma^*$ leads to one last condition for the $CSI$ spacetime.  Since $\ov{W}{1}{e} = -(ln \zeta)_{,e}$, the only constraint on the function $\zeta$ arises from \eqref{Wcsi4}.

\beq \displaystyle \sum_{i=3}^{N} [D_i ln(\zeta_2) ]^2 = -4\sigma,~~~~ \sigma < 0 \label{c112zetaconstraint} \eeq

\beq g^{ef}\ov{W}{1}{e} \ov{W}{1}{f} = - 4 \sigma. \nonumber \eeq

\noindent  The left-hand-side must be positive, and so it is necessary that $\sigma = R_{1212}$ is a negative real number. 

\subsection*{Case 1.2 :}

\noindent The remaining conditions from $|\zeta| \leq 0$ are 

\beq &D_2 \zeta_1 - \ov{H}{1}{~} \zeta_1 = \ov{W}{1}{3} \zeta_3& \label{c12deqn} \\
& (\zeta_3)^2 \leq \zeta_1 \ov{\zeta}{0}{2}&. \label{c12inneq} \eeq

\noindent  Expanding $\zeta_2$ and $\ov{\Gamma}{1}{3}$, we find that the $O(v^2)$ terms, \eqref{2ke23} and \eqref{2ke2n} are automatically satisfied, while \eqref{2ke22} and using \eqref{c12deqn} yield the following differential equation for $\sigma^*(u,x^e)$:

\beq \zeta_1  D_2 \sigma^* + \zeta_3 D_3 \sigma^*  = 0. \label{c12sigmaeqn} \eeq

Using a coordinate transformation of the form \eqref{type3ct}, coordinates are chosen  so that $\zeta_1 = 1$ and the non-spacelike condition \eqref{c12deqn} determines a part of $H$

\beq \ov{H}{1}{~} = - \ov{W}{1}{3} \zeta_3. \nonumber \eeq

\noindent We can apply another coordinate transform of type \eqref{type2ct} to eliminate $\ov{H}{0}{~}$ as well. In this coordinate system the Killing equations are:

\beq & \sigma^* ( \ov{\zeta}{0}{2} - \zeta_3 \ov{W}{0}{3} ) = 0, & \label{keqn2}   \\
 & D_2 \ov{\zeta}{0}{2} + \zeta_3 D_3 \ov{\zeta}{0}{2} + \zeta_3 D_2\zeta_3 = 0, & \label{keqn3} \\
 & D_2 \ov{W}{1}{i} + \zeta_3 D_3 \ov{W}{1}{i} = 0, & \label{keqn4} \\
 & D_2 \ov{W}{0}{3} - \zeta_3 \ov{W}{0}{3} \ov{W}{1}{3} = -D_2 \zeta_3 - D_3 \ov{\zeta}{0}{2} - \ov{W}{1}{3} \ov{\zeta}{0}{2},  & \label{keqn5} \\
 & D_2 \ov{W}{0}{n} + \zeta_3 D_3 \ov{W}{0}{n} = (\ov{W}{0}{3} \ov{W}{1}{n} + D_n \ov{W}{0}{3}) \zeta_3 - \ov{W}{1}{n} \ov{\zeta}{0}{2} - D_n \ov{\zeta}{0}{2}, & \label{keqn6} \eeq

\noindent If $\zeta_3$ is non-zero, equation \eqref{keqn3} simplifies the differential equation for $\ov{W}{0}{3}$ in \eqref{keqn5}. Thus two subcases must be considered in which $\zeta$ vanishes or not.

\subsection*{Case 1.2.1: }

Setting $\zeta_3$ equal to zero causes $\ov{H}{1}{~}$ to vanish while \eqref{keqn3} and \eqref{keqn4} imply

\beq & D_2 \ov{W}{1}{i} = D_2 \ov{\zeta}{0}{2} = 0.& \label{c121keqn1} \eeq 

\noindent The remaining equations give constraints for the remaining metric functions:

\beq & \sigma^* ( \ov{\zeta}{0}{2} ) = 0, & \label{c121keqn2}   \\
 & D_2 \ov{W}{0}{3} = - D_3 \ov{\zeta}{0}{2} - \ov{W}{1}{3} \ov{\zeta}{0}{2},  & \label{c121keqn3} \\
 & D_2 \ov{W}{0}{n} = - \ov{W}{1}{n} \ov{\zeta}{0}{2} - D_n \ov{\zeta}{0}{2}. & \label{c121keqn4} \eeq

\noindent Thus there are two minor subcases to consider arising from \eqref{c121keqn2}.

\subsection*{ Case 1.2.1a } Assuming $\sigma^* \neq 0$,  $\ov{\zeta}{0}{2}$ vanishes and the set of spacetimes with metric functions: 

\beq H(v,x^e) = \sigma^* \frac{v^2}{8},~~~W_i(v,x^e) = \ov{W}{1}{i}(x^e)v + \ov{W}{0}{i}(x^e) \label{c121metricfnA} \eeq%

\noindent are $CCNV$ spacetimes with $\frac{\partial}{\partial u}$ as a covariantly constant null vector admitting a Killing vector of the form:

\beq \zeta = n + \frac{\sigma^* v^2}{8} \ell. \nonumber \eeq



If we suppose $\zeta$ is a covariantly constant vector; Lemma \eqref{lem:azeta0} and equations \eqref{azeta23} and \eqref{azeta2n} force the metric functions $\ov{W}{1}{i}$ and $\ov{W}{0}{n}$ to vanish. However a contradiction arises from \eqref{azeta12} as it requires $\sigma^* = 0$ but we have assumed that $\sigma^* \neq 0$ and so the above spacetime cannot admit a covariantly constant vector.

\subsection*{Case 1.2.1b} For the other subcase, $\sigma^*$ is equal to zero, and the positive-definite signature of the transverse metric restricts $\sigma \leq 0$. For arbitrary $\ov{\zeta}{0}{2}(x^e)$, and any choice of $\ov{W}{1}{i}(x^e)$ satisfying \eqref{c11metric} with $\sigma = R_{1212} \leq 0$, the $CSI$ Kundt spacetime with a locally homogeneous transverse space and metric functions:

\beq H = 0,~~~W_i(u,v,x^e) = \ov{W}{1}{i}v - (D_i\ov{\zeta}{0}{2} + \ov{W}{1}{i}\ov{\zeta}{0}{2})u + w_i(x^e)  \label{c121metricfnB} \eeq

\noindent admit a Killing vector of the form:

\beq \zeta = n + \ov{\zeta}{0}{2} \ell \nonumber \eeq

\noindent To preserve the non-spacelike requirement $\ov{\zeta}{0}{2}$ must always be greater than or equal to zero. If this killing vector is covariantly constant, $\ov{W}{1}{i} = 0$ and hence $\sigma = 0$, equation \eqref{azeta3n} implies $A_{ij} = 0$, and the remaining equations \eqref{azeta23} and \eqref{azeta2n} force $\ov{\zeta}{0}{2}$ to be constant. Thus $\zeta$ is the sum of the $CCNV$'s $\ell$ and $n$.

\subsection*{Case 1.2.2: $\zeta_3 \neq 0$} 

\noindent Divide by $\zeta_3$ in \eqref{keqn3} and substitute the result into \eqref{keqn5} to simplify the differential equation for $\ov{W}{0}{3}$:

\beq D_2 \ov{W}{0}{3} - \zeta_3 \ov{W}{0}{3} \ov{W}{1}{3} = \frac{D_2 \ov{\zeta}{0}{2} }{\zeta_3} - \ov{W}{1}{3} \ov{\zeta}{0}{2} \label{c122w03A} \eeq

\noindent then by multiplying the above by $E(u,x^e) = e^{- [ \int \ov{W}{1}{3} \zeta_3 du] }$, integration by parts gives the solution

\beq \ov{W}{0}{3} = \frac{ \ov{\zeta}{0}{2} }{\zeta_3} + e^{ [ \int \ov{W}{1}{3} \zeta_3 du]} \int \frac{ \ov{\zeta}{0}{2} D_2 \zeta_3   }{(\zeta_3)^2 e^{ [ \int \ov{W}{1}{3} \zeta_3 du]} } du. \nonumber \eeq

\noindent From \eqref{keqn2} there are two minor subcases to consider, depending upon whether $\sigma^*$ vanishes or not. 

\subsection*{Case 1.2.2a :} Supposing that $\sigma^*$ does indeed vanish, the functions $\ov{W}{1}{3}(x^e)$ and $\ov{W}{1}{n}$ must satisfy \eqref{c11metric} with $\sigma \leq 0$. For arbitrary $\zeta_3(u)$ and any solution of the following differential equation 

\beq  & D_2 \ov{\zeta}{0}{2} + \zeta_3 D_3 \ov{\zeta}{0}{2} = - \zeta_3 D_2\zeta_3  & \label{c122zeta2A} \eeq

\noindent the Kundt $CSI$ spacetime with a locally homogeneous transverse space and 

\beq & H = - \ov{W}{1}{3} \zeta_3 v &  \nonumber \\
 &  W_3(u,v,x^e) = \ov{W}{1}{3}(u,x^e) v + \frac{ \ov{\zeta}{0}{2} }{\zeta_3} + \frac{1}{E}\int \frac{ E \ov{\zeta}{0}{2} D_2 \zeta_3 }{(\zeta_3)^2 } du,~~~~E = e^{\int \ov{H}{1}{~} du} &   \label{c122metricfnA} \\
&   W_n(u,v,x^e) = \ov{W}{1}{n}(u,x^e)v + \ov{W}{0}{n}(u,x^e) &  \nonumber \eeq 
 
\noindent satisfying the following differential equations:

\beq & D_2 \ov{W}{1}{i} + \zeta_3 D_3 \ov{W}{1}{i} = 0 & \label{c122w1nA} \\
 & D_2 \ov{W}{0}{n} + \zeta_3 D_3 \ov{W}{0}{n} = \frac{\zeta_3 \ov{W}{1}{n}}{E}\int \frac{ E \ov{\zeta}{0}{2} }{(\zeta_3)^2 } du + D_n[ \frac{\zeta_3}{E}\int \frac{ E \ov{\zeta}{0}{2} D_2 \zeta_3 }{(\zeta_3)^2 } du ] \label{c122w0nA} \eeq

\noindent admits a Killing vector of the form

\beq \ell + \ov{\zeta}{0}{2}(u,x^e) n + \zeta_3(u) m^3 \nonumber \eeq

Requiring  $\zeta$ to be a $CCNV$, the $\ov{W}{1}{i}$ must vanish, causing $H = 0$ and $\sigma = 0$.  This is an example of a $CCNV$ metric, with $\ell = \frac{\partial}{\partial{v}}$ as the $CCNV$, where $\zeta$ will be a second $CCNV$. The additional constraints \eqref{azeta12} - \eqref{azeta3n} imply $A_{ij} = 0$ while the remaining equations lead to two possible subcases for Kundt spacetimes admitting a covariantly constant vector, either $D_n\zeta_2 = 0$ or $D_2 \zeta_3 = 0$. The first case leads to the following form for $\zeta$ and the metric functions

\beq & \zeta = n + [-\zeta^2_3] \ell + \zeta_3(u) m^3,~~~~ 3\zeta_3^2 \leq 0 & \nonumber\\
& H = 0,~~~~ W_3(u,x^e) = -\zeta_3 + w_3(x^e),~~~~  W_n(x^e) = \int D_n w_3 dx^3 + w_n( x^r).& \nonumber \eeq 


\noindent The non-spacelike condition $3 \zeta_3^2 \leq 0$ eliminates the above case, as we've assumed $\zeta_3 \neq 0$ this case is not admissible. In the second case $\zeta_3$ must be constant, scaling $x^3$ so that $\zeta_3 = 1$,

\beq & \zeta = n + \zeta_2(x^r) \ell + m^3,~~~~ 1 \leq 2 \zeta & \nonumber\\
& H = 0,~~~~ W_3(x^e) =  w_3(x^e),~~~~  W_n(x^e) = \int D_n (w_3) dx^3 - 2 D_n ( \zeta_2 )x^3 + w_n( x^r),& \nonumber \eeq 

\noindent The vanishing of $A_{3j} =  D_{[j}\ov{W}{0}{3]}$ implies that $D_n\ov{\zeta}{0}{2} =0$ and so $\ov{\zeta}{0}{2}$ must be constant. If $\zeta$ is timelike, the constant  $\ov{\zeta}{0}{2} > \frac12$ while if $\zeta$ is null $\ov{\zeta}{0}{2} = \frac12$. These spacetimes will automatically be $CCNV$ spacetimes with $\ell$ as another covariantly constant null vector. 

\subsection{Case 1.2.2b:} If $\sigma^* $ is non-zero, it satisfies the differential equation \eqref{c12sigmaeqn}

\beq  D_2 \sigma^* + \zeta_3 D_3 \sigma^*  = 0 \nonumber \eeq 

\noindent and the identity  $\ov{\zeta}{0}{2} = \zeta_3 \ov{W}{0}{3}$ may be derived from \eqref{keqn2}, which causes \eqref{c122w03A} to simplify, implying $ \ov{\zeta}{0}{2} D_2 \zeta_3 = 0$. Letting $\ov{\zeta}{0}{2} = 0$, the differential equation \eqref{c122zeta2A} for $\ov{\zeta}{0}{2}$ forces $D_2 \zeta_3 = 0$. In either case $\zeta_3$ must be constant and henceforth will be set to one. For any solution $\ov{\zeta}{0}{2}$ to the differential equation

\beq D_2\ov{\zeta}{0}{2} + D_3\ov{\zeta}{0}{2} = 0 \label{c122zeta2B}, \eeq

\noindent the vector  

\beq n + [\frac{\sigma^*}{8}v^2 + \ov{\zeta}{0}{2}] \ell + m^3 \nonumber \eeq

\noindent will be a Killing vector for any $CSI$ Kundt spacetime of the form

\beq & H = \frac{\sigma^*}{8}v^2 - \ov{W}{1}{3} \zeta_3 v, & \label{c122metricfnB} \\
& W_3(u,v,x^e) = \ov{W}{1}{3}(u,x^e) v + \ov{\zeta}{0}{2},~~ W_n(u,v,x^e) = \ov{W}{1}{n}(u,x^e)v + \ov{W}{0}{n}(u,x^e) & \nonumber \eeq 
 

\noindent where the $\ov{W}{1}{i}$ and $\ov{W}{0}{i}$ satisfy the following equations:

\beq D_2 \ov{W}{1}{i} + D_3 \ov{W}{1}{i} & = & 0,   \label{c122w1nB} \\
  D_2 \ov{W}{0}{i} + D_3 \ov{W}{0}{i} & = & 0.  \label{c122w0nB} \eeq

\noindent If $\zeta$ is required to be covariantly constant, a contradiction arises from \eqref{azeta12} as it requires $\sigma^* = 0$ despite the fact that we have assumed $\sigma^* \neq 0$. Thus there are no $CCNV$ spacetimes of the form \eqref{c122metricfnB}.
 
\section*{Case 2 : $\ov{\Gamma}{2}{3} = 0$ } 







For the remainder of this case we shall assume $D_3\zeta_1 \neq 0$ to avoid the previous subcases.  Supposing $\ov{W}{1}{3} = 0$, this implies that $D_3D_3 \zeta_1 =0$ and $\sigma^* =\sigma$.  This causes a contradiction to arise between the  Killing equation \eqref{2ke23} and the non-spacelike condition \eqref{magv2}: 

\beq & 2\sigma D_3(\zeta_1) = 0,~~~(D_3\zeta_1)^2 \leq \sigma (\zeta_1)^2.& \nonumber\eeq

\noindent The first implies that $\sigma = 0$ as we have assumed $D_3\zeta_1 \neq 0$; however, by the second inequality the vanishing of $\sigma$ implies $D_3\zeta_1 = 0$ which contradicts our original assumption. Thus $D_3D_3 \zeta_1$ is always non-zero, and using this fact we may derive another identity for  $\sigma^*= 4\sigma + (\ov{W}{1}{3})^2$ in terms of $\zeta_1$ from the vanishing of $\ov{\Gamma}{2}{3}$: 

\beq \sigma^* = \frac{D_3 \sigma^*}{\ov{W}{1}{3}} = \frac{ 2 \ov{W}{1}{3} D_3\ov{W}{1}{3} }{ \ov{W}{1}{3} } = 2D_3D_3(lnD_3 \zeta_1). \label{c2sigmastar} \eeq

\noindent Using a coordinate transform of type \eqref{type3ct} with $g(u) = \frac{u}{\sqrt{|\sigma|}}$, we may rescale $\sigma$ in  \eqref{HWiSigma} so that it equals $\sigma = -1,0,1$ depending on it's sign. Doing so will scale all of the metric functions and Killing vector components by a constant value, but otherwise will leave them unchanged.  

Dropping the primes and substituting \eqref{c2sigmastar} into the original identity for $\sigma^*$ yields another differential equation for $D_3 \zeta_1$:

\beq  D_3D_3ln (D_3 \zeta_1) - \frac12 ( D_3ln( D_3 \zeta_1 ) )^2 = 2 \sigma . \nonumber \eeq

\noindent Multiplication by $exp( - \frac12 \int D_3(ln D_3 \zeta_1) dx^3) = (D_3 \zeta_1)^{-\frac12}$ leads to the simpler equation

\beq  D_3D_3[ (D_3\zeta_1)^{-\frac12}] =  - \sigma (D_3\zeta_1)^{-\frac12}. \label{gamma23eqn} \eeq

\noindent There are three possible solutions to this equation depending on whether $\sigma$ is positive, negative or zero: 

\beq \sigma = -1 &:& (D_3\zeta_1)^{-\frac12} = c_1(u) cosh(x^3) + c_2(u) sinh(x^3),  \nonumber \\ 
\sigma = 0 &:& (D_3\zeta_1)^{-\frac12} = c'_1(u)x^3 + c'_2,  \nonumber \\
\sigma = 1 &:& (D_3\zeta_1)^{-\frac12} = c''_1(u) cos(x^3) + c''_2(u) sin(x^3).  \nonumber \eeq

Ignoring these facts for a moment, we recall that the metric functions $W_i$ may be expressed in terms of $\ov{\zeta}{0}{3}$ and $\zeta_1$ using \eqref{zeta3iskey} :

\beq \ov{W}{0}{i} = \frac{ - D_i \ov{\zeta}{0}{3} }{ D_3 \zeta_1 }, ~~~ \ov{W}{1}{i} = D_i ln ( D_3 \zeta_1).  \nonumber \eeq

\noindent In this case, it is possible to set all but $\ov{W}{1}{3}$ to zero by making a coordinate transform of type \eqref{type2ct} with $h = - \frac{\ov{\zeta}{0}{3}}{D_3\zeta_1}$. In these new coordinates, the metric functions take the form:

\beq &W_3 = D_3ln(D_3\zeta_1) v,~~ W_n = 0.& \label{c2metricfn} \eeq

\noindent The following coefficient functions of $H$ change in the new coordinate system:

\beq & \ov{H}{1}{~} = \ov{H'}{1}{~} + \frac{\ov{\zeta}{0}{3} \sigma^*}{4 D_3\zeta_1},~~ \ov{H}{0}{~} = \ov{H'}{0}{~} +  \frac{\ov{\zeta}{0}{3} \ov{H'}{1}{~}}{ D_3\zeta_1}  + D_2 \left(\frac{\ov{\zeta}{0}{3} }{ D_3\zeta_1} \right)  \frac{\sigma^* (\ov{\zeta}{0}{3})^2}{8 (D_3\zeta_1)^2}, \nonumber \eeq

\beq & \ov{H}{1}{~} = \ov{H'}{1}{~} + \frac{\ov{\zeta}{0}{3} \sigma^*}{4 D_3\zeta_1},~~ \ov{H}{0}{~} = \ov{H'}{0}{~} +  \frac{\ov{\zeta}{0}{3} \ov{H'}{1}{~}}{ D_3\zeta_1}  + D_2 \left(\frac{\ov{\zeta}{0}{3} }{ D_3\zeta_1} \right)  \frac{\sigma^* (\ov{\zeta}{0}{3})^2}{8 (D_3\zeta_1)^2}, \nonumber \eeq 

\noindent where primed functions denote the functions in the previous coordinate system. As the original $\ov{H'}{1}{~}$ and $\ov{H'}{0}{~}$ were arbitrary functions of $u$ and the spatial coordinates, we may ignore the special form the v-coefficients take in this coordinate system and treat them simply as new arbitrary functions. In this coordinate system the tensor $A_{3n}$ given in \eqref{note1} vanishes, and the connection coefficients $\Gamma_{2i2}$ are of the form:

\beq \Gamma_{2i2} &=& \overbrace{(D_i \ov{H}{1}{~} - D_2\ov{W}{1}{i})}^{ \ov{\Gamma}{1}{i} } v  + \overbrace{ D_i \ov{H}{0}{~} + \ov{H}{0}{~} \ov{W}{1}{i} }^{ \ov{\Gamma}{0}{i} } \nonumber. \eeq

\noindent This choice of coordinate system simplifies the Killing equations considerably; for example, the other two covector components are now

\beq & \zeta_2 = \overbrace{ ( \frac{ \sigma^* \zeta_1}{4} - D_3D_3\zeta_1 ) }^{\ov{\zeta}{2}{2} } \frac{v^2}{2} + \overbrace{ ( \ov{H}{1}{~} \zeta_1  - D_2 \zeta_1 ) }^{\ov{\zeta}{1}{2} } v + \ov{\zeta}{0}{2}(u,x^e), & \nonumber \\
& \zeta_3 = - D_3 (\zeta_1) v. \nonumber \eeq

\noindent Taking the magnitude of the vector and invoking the non-spacelike conditions yield 

\beq & D_3D_3ln[ (D_3 \zeta_1)^{-\frac12}] + D_3(ln(D_3\zeta_1))D_3ln(\zeta_1) + (D_3ln(\zeta_1))^2 \leq 0, & \label{c2magv2}\\
&  \zeta_1 (\ov{H}{1}{~} \zeta_1 - D_2 \zeta_1 ) = 0, & \label{c2magv1}\\
& \zeta_1 \ov{\zeta}{0}{2} \geq 0. & \label{c2magv0} \eeq

\noindent Thus $\ov{\zeta}{1}{2}$ must vanish and we may solve for $\ov{H}{1}{~}$ in terms of $\zeta_1$,

\beq \ov{H}{1}{~} = D_2 ln(\zeta_1). \nonumber \eeq

\noindent Further constraints on $H$ involving $\ov{H}{0}{~}$ may be found by taking those Killing equations involving the spatial derivatives of $\ov{\zeta}{0}{2}$; i.e., \eqref{0ke23} and \eqref{0ke2n} and considering integrability conditions. In this coordinate system \eqref{0ke23} and \eqref{0ke2n} are 

\beq & D_3 \ov{\zeta}{0}{2} + \ov{H}{0}{~} D_3\zeta_1 - \zeta_1 D_3 \ov{H}{0}{~} - \zeta_1 \ov{H}{0}{~} D_3ln(D_3\zeta_1) + \ov{\zeta}{0}{2} D_3ln(D_3\zeta_1) = 0 & \nonumber \\
&D_n \ov{\zeta}{0}{2} - \zeta_1 D_n \ov{H}{0}{~} = 0 & \nonumber \eeq

\noindent We note that the commutator applied to any function independent of $v$ vanishes (i.e. $[D_3, D_n] f(u,x^e) = 0$); thus   differentiating the first equation by $D_n$ and the latter by $D_3$ and subtracting the result gives the following constraint

\beq 2 D_n(\ov{H}{0}{~}) D_3\zeta_1  = 0. \nonumber \eeq

\noindent Hence $\ov{H}{0}{~}$ and $\ov{\zeta}{0}{2}$ are actually  functions of $u$ and the spatial coordinate $x^3$. 

In light of this fact the Killing equations \eqref{2ke2n} - \eqref{0ke2n} are automatically satisfied. Similarly,  equation \eqref{2ke23} may be ignored as it gives the identity $\sigma^* = 2D_3D_3ln(D_3\zeta_1)$, which arose from the vanishing of $\ov{\Gamma}{2}{3}$. The remaining Killing equations are now:

\beq &D_2 \sigma^* = 4 D_2 \left( \frac{D_3D_3 \zeta_1}{\zeta_1} \right) - \frac12 D_2[ (D_3ln(\zeta_1))^2 ],  & \label{c2-2ke22} \\
& D_3 \ov{H}{0}{~} = \frac{\sigma^*}{4D_3\zeta_1}(\ov{\zeta}{0}{2} - \ov{H}{0}{~} \zeta_1), & \label{c2-1ke22} \\
& D_2 \ov{\zeta}{0}{2} = -\ov{\zeta}{0}{2} D_2ln(\zeta_1), & \label{c2-0ke22} \\
& 2D_2D_3ln(\zeta_1) = D_2D_3 ln(D_3\zeta_1), & \label{c2-1ke23} \\
& D_3( \ov{\zeta}{0}{2} D_3\zeta_1) = \zeta_1^2 D_3[ \ov{H}{0}{~} D_3ln(\zeta_1)]. & \label{c2-0ke23} \eeq

\noindent Differentiating \eqref{c2-1ke23} and using the fact that $[D_3, D_2] f(u,x^e) = 0$, one finds the following expression for $D_2 \sigma^* = 2D_2D_3D_3ln(D_3\zeta_1)$:

\beq D_2 \sigma^* = 4 D_2 \left[ \frac{D_3D_3 \zeta_1}{\zeta_1} - \left( \frac{D_3 \zeta_1}{\zeta_1} \right)^2 \right]. \nonumber \eeq

\noindent Subtracting this from \eqref{c2-2ke22} yields the following constraint

\beq D_2 \left( D_3ln(\zeta_1) \right)^2 = 0 \nonumber \eeq

\noindent implying that $\zeta_1$ must take the form:

\beq \zeta_1 = e^{A(x^3)} e^{B(u)} \label{c2-zeta1}. \eeq

\noindent Apply a coordinate transform of type \eqref{type3ct} with $g = \int e^{-B(u)} du$ will remove the $u$ dependence from $\zeta_1$. Rewriting \eqref{c2-0ke22} in terms of $\ov{\zeta'}{0}{2} = \ov{\zeta}{0}{2}e^B$, it is easily shown that this implies  $D_2 \ov{\zeta'}{0}{2} = 0$. Denoting $\zeta'_1 = e^{A(x^3)}$ the Killing vector $\zeta =  e^A e^B n + \zeta_2 \ell + \zeta_3 m^3$ becomes:

\beq \zeta = \zeta'_1 n' + \left[ \left( \frac{\sigma^* \zeta'_1}{4} - D_3D_3(\zeta'_1) \right) \frac{{v'}^2}{2} + \ov{\zeta'}{0}{2}(x^3) \right] \ell' + [-D_3(\zeta'_1) v'] m^3 \label{c2kvform} . \nonumber \eeq

\noindent In the remaining Killing equations,  \eqref{c2-1ke22} and \eqref{c2-0ke23}, the function $\ov{H}{0}{~}$ in the new coordinate system becomes  $\ov{H'}{0}{~} = e^{2B} \ov{H}{0}{~}$ and so we may remove $e^B$ entirely from these two equations. 

Dropping the primes and combining \eqref{c2-1ke22} with \eqref{c2-0ke23} yields the following algebraic equation for $\ov{H}{0}{~}$: 

\beq \ov{H}{0}{~} \left( D_3D_3ln(\zeta_1) - \frac{\sigma^*}{4} \right) = \frac{D_3( D_3(\zeta_1) \ov{\zeta}{0}{2} )}{\zeta_1^2} - \frac{\sigma^* \ov{\zeta}{0}{2} }{4 \zeta_1} \nonumber \eeq

\noindent The coefficient of $\ov{H}{0}{~}$ cannot vanish, as the non-spacelike condition \eqref{c2magv2} would imply 

\beq & 2(D_3 ln(\zeta_1))^2 \leq 0. & \nonumber \eeq

\noindent It is assumed that $D_3\zeta_1 \neq 0$ so the above constraint is impossible. Simplifying the above expression $\ov{H}{0}{~}$ may be written as

\beq  \ov{H}{0}{~} = \frac{D_3( D_3(\zeta_1) \ov{\zeta}{0}{2} ) + D_3D_3ln((D_3\zeta_1)^{-\frac12}) \ov{\zeta}{0}{2} \zeta_1 }{ \zeta_1^2 D_3D_3 ln( \zeta_1 (D_3 \zeta_1)^{-\frac12})} \label{c2-H0} \eeq

\noindent Having exhausted the Killing equations, we look to the remaining non-spacelike conditions \eqref{c2magv2} and \eqref{c2magv0}.
\subsection*{ Case 2.1: Null Killing Vectors}

\noindent If $\zeta$ is required to be null $\ov{\zeta}{0}{2}$ must be zero, forcing $\ov{H}{0}{~}$ to vanish as well. Using \eqref{c2magv2} and \eqref{gamma23eqn} we find the following expression

\beq & D_3(A) = D_3ln[ (D_3 \zeta_1)^{-\frac12} ] \pm \sqrt{  2[ D_3ln( (D_3 \zeta_1)^{-\frac12} )]^2 + \sigma}. & \label{c2nulleqn} \eeq

\noindent Combining this with the solution to \eqref{gamma23eqn} for a particular $\sigma =-1,0,1$:

\beq \sigma = -1 &:& (D_3\zeta_1)^{-\frac12} = c_1 cosh(x^3) + c_2 sinh(x^3)  \label{c2d3zeta1neg} \\ 
\sigma = 0 &:& (D_3\zeta_1)^{-\frac12} = c_1x^3  \label{c2d3zeta1zero}+ c_2  \\
\sigma = 1 &:& (D_3\zeta_1)^{-\frac12} = c_1 cos(x^3) + c_2 sin(x^3)  \label{c2d3zeta1pos} \eeq

\noindent we may algebraically solve for $\zeta_1$ by noting that $D_3 \zeta_1 = D_3(A) e^A = D_3(A) \zeta_1$: 



\beq & \sigma = -1 : \zeta_1 = \frac{ (c_1 cosh(x^3) + c_2 sinh(x^3))^{-1} }{  c_1 sinh(x^3) + c_2 cosh(x^3) \pm \sqrt{ c_1^2 + c_2^2 + \left(c_1 sinh(x^3) + c_2cosh(x^3) \right)^2 } } & \label{c2nullneg} \\
\nonumber \\
& \sigma = 0 : \zeta_1 = \frac{ 1}{ c_1 (1 \pm \sqrt{2}) (c_1 x^3 + c_2)} & \label{c2nullzero} \\
\nonumber \\
& \sigma = 1 : \zeta_1 = \frac{(c_1 cos(x^3) + c_2 sin(x^3))^{-1}}{  -c_1 sin(x^3) + c_2 cos(x^3) \pm \sqrt{ c_1^2 + c_2^2 + \left(- c_1 sin(x^3) + c_2 cos(x^3) \right)^2 } } & \label{c2nullpos} \eeq

Supposing that $\zeta$ is covariantly constant, the constraint in Lemma \eqref{lem:azeta0} on $\ov{W}{1}{3}$ along with the identity \eqref{zeta3iskey} yields

\beq D_3ln(\zeta_1) = - D_3ln[(D_3 \zeta_1)^{\frac12}]. \label{ccnvcondition} \eeq

\noindent Since  $ln(\zeta_1) = A$, the above simplifies \eqref{c2nulleqn} in the null case, giving

\beq  2 D_3ln[ (D_3 \zeta_1)^{-\frac12} ] \pm \sqrt{  2[ D_3ln( (D_3 \zeta_1)^{-\frac12} )]^2 + \sigma} = 0 \nonumber \eeq

\noindent Multiplying both roots together the result must vanish

\beq 2[D_3ln[ (D_3 \zeta_1)^{-\frac12} ]]^2 - \sigma = 0 \label{c2ccnv} \eeq 

\noindent Substituting the three posibilities of $(\zeta_1)^{\frac12}$ gives the constraint: 

\beq & 3[c_1^2 - c_2^2]sinh^2(x^3) +6c_1c_2sinh(x^3)cosh(x^3)+2c_2^2 + c_1^2 = 0 & \nonumber \\ 
 & \frac{c_1^2}{(c_1 x^3 + c_2)^2} = 0 & \nonumber \\
 & \sigma = 1 : 3[c_1^2 - c^2] sin^2(x^3) -6c_1c_2sin(x^3)cos(x^3) + 2c_2^2 - c_1^2 = 0 & \nonumber \eeq

\noindent In each case this identity will only hold if $c_1 = c_2 = 0$; however, this will imply that $D_3\zeta_1 = 0$, which cannot happen. Thus the null killing vector $\zeta$ cannot be covariantly constant.

\subsection*{ Case 2.2: Timelike Killing Vectors}

\noindent If we require $\zeta$ to be timelike, equation \eqref{c2magv0} along with the fact that $\zeta_1 = e^A$ forces $\ov{\zeta}{0}{2}$ to be greater than or equal to zero for all values of $x^3$. To find $\zeta_1$ we integrate each of the three solutions to \eqref{gamma23eqn} given above

\beq & \sigma = -1 : \zeta_1 = \frac{sinh(x^3)}{c_1(c_1cosh(x^3) + c_2sinh(x^3))} + c_3 & \label{c2tlneg} \\
\nonumber \\
& \sigma = 0 : \zeta_1 = \frac{-1}{c_1(c_1 x^3 + c_2)} + c_3 \label{c2tlzero} \\
\nonumber \\
& \sigma = 1 : \zeta_1 = \frac{sin(x^3)}{c_1(c_1cos(x^3) + c_2sin(x^3))} + c_3. & \label{c2tlpos} \eeq

\noindent The inequality \eqref{c2magv2} restricts the choice of $c_3$  depending on the choice of $c_1$ and $c_2$: 

\beq & \sigma = -1 : [c_1^2 + c_2^2] \zeta_1^2 - 2 \left( \frac{ c_1 sinh(x^3) + c_2 cosh(x^3)}{c_1 cosh(x^3) + c_2 sinh(x^3) } \right) \zeta_1 + \frac{1}{(c_1 cosh(x^3) + c_2 sinh(x^3))^2} < 0 & \nonumber \\
\nonumber \\
& \sigma = 0 : -c_1^2 \zeta_1^2 - 2 \left( \frac{c_1}{c_1 x^3 + c_2} \right) \zeta_1 + \frac{1}{(c_1x^3 + c_2)^2} < 0 & \label{c2tlinneq} \\
\nonumber \\
& \sigma = 1 : -[c_1^2 + c_2^2] \zeta_1^2 - 2 \left( \frac{ -c_1 sin(x^3) + c_2 cos(x^3) }{c_1 cos(x^3) + c_2 sin(x^3)} \right)  \zeta_1 + \frac{1}{(c_1 cos(x^3) + c_2 sin(x^3))^2} < 0 & \nonumber \eeq


\noindent Notice in both the null and timelike case, the value of $\sigma$ restricts the domain of $x^3$. When $\sigma = 1$, the domain of $x^3$ is limited to a finite interval, $x^3 \in (x^3_0, x^3_0+ \pi)$, as the value $x^3_0 = arctan(- \frac{c_1}{c_2}) $ will cause $(D_3 \zeta_1)^{-\frac12}$ to vanish. When $\sigma = 0$, $x^3 \geq -\frac{c_2}{c_1}$ to avoid singularities. In the case with $\sigma = -1$, $x^3 > x^3_0 = arctanh(-\frac{c_1}{c_2})$ when $c_1/c_2 \leq 1 $, otherwise $\zeta_1$ is regular on the whole of the real line.

\noindent Requiring $\zeta$ to be covariantly constant, equation \eqref{ccnvcondition} may be rewritten as a function set to zero in terms of $\zeta$ and $(D_3\zeta_1)^{\frac12}$ for the three subcases with $\sigma = -1,0,1$ respectively:  

\beq & [c_1 + 2c_1^2 c_2 c_3]sinh^2(x^3) + [c_2 + c_1^3 c_3 + c_1 c_2^2 c_3]cosh(x^3)sinh(x^3)+ [c_1^2 c_2 c_3 + c_1] & \nonumber \\
&  c_1 c_3 (c_1x^3 + c_2) & \nonumber \\
&  -[c_1 + 2 c_1^2 c_2 c_3]sin^2(x^3) + [c_2 - c_1^3 c_3 + c_1 c_2^2 c_3]cos(x^3)sin(x^3) + [c_1^2 c_2 c_3 + c_1]  & \nonumber \eeq

\noindent In both cases where $\sigma = -1,1$ the vanishing of the first and third equation  will hold only if $c_1$ and $c_2$ both vanish, which violates the assumption $D_3 \zeta_1 \neq  0$, and so there are no timelike covariantly constant vectors in either of these two cases. When $\sigma = 0$, setting the second equation to zero implies $c_3 =0$, the Killing vector of the form \eqref{c2kvform} with $\zeta_1 = -1/(c_1^2 x^3 + c_1 c_2)$ satsifies the condition in \eqref{ccnvcondition}. A problem arises from the inequality \eqref{c2magv2}

\beq & -c_1^2 \left( \frac{1}{c_1^2(c_1x^3 + c_2)^2} \right) - 2 \left( \frac{c_1}{ (c_1 x^3 + c_2)} \right) \left( \frac{-1}{c_1(c_1x^3 + c_2)} \right) + \frac{1}{(c_1 x^3 + c_2)^2} < 0;&  \nonumber \eeq

\noindent simplifying the above leads to the inequality $2 < 0$ which is clearly impossible. We conclude there are no covariantly constant timelike vectors in the spacetimes belonging to Case 2. 
 
\section*{Conclusions}

To determine the subset of Kundt $CSI$ spacetimes admitting a null or timelike isometry, several choices were made to simplify the Killing equations. Local coordinates were chosen so that one of the spacelike Killing vectors, $Y$, belonging to the (locally) homogeneous transverse space has been rectified to act locally as a translation in the $x^3$ direction,  i.e., $Y = A \frac{\partial}{\partial x^3}$. The frame was then rotated so that the frame vector $m^3$ was aligned with the spatial part of $\zeta$ and, moreover, that the matrix $m_{ie}$ was upper-triangular with $m_{33} = 1$. This causes the connection components $\Gamma_{3ij}$ and $\Gamma_{ij3}$ to vanish, simplifying the Killing equations considerably. 

In this coordinate system we determined the special form for the components of $\zeta$ in terms of arbitrary functions and in terms of $H$ and the $W_e$;  i.e., \eqref{kvcomp13} and \eqref{kvcomp2}. All of the functions involved (metric or otherwise) are expressed as polynomials in $v$ with coefficient functions of $u$ and $x^e$. These are substituted into the remaining Killing equations which are rearranged into the various orders of $v$ to give \eqref{3ke22} - \eqref{0ke2n}, while the non-spacelike conditions yield \eqref{magv2} - \eqref{magv0}. The highest order equation \eqref{3ke22} gives two major subcases, either $D_3 \zeta_1 = 0$ or $\ov{\Gamma}{2}{3}=0$ in \eqref{note1}. 

It is known that  all $VSI$ spacetimes admitting a non-spacelike isometry are $CCNV$ spacetimes with $\ell$ as the covariantly constant vector \cite{VSISugra}. As an analogue to this result, the equations arising from $\nabla_{[a} \zeta_{b]} = 0$ were examined to determine which $CSI$ Kundt spacetimes admit a covariantly constant vector and which cannot. 

The results of the analysis are summarized below:

\subsection*{ Case 1.1.1: $\zeta = \ell$}

\noindent In this case $R_{1212} = \sigma = 0$, the metric functions in \eqref{Kundt} takes the form: $H(u,x^k)$ and $W_{ i}(u,x^k)$. All $CSI$ spacetimes in this subcase are clearly $CCNV$ spacetimes with $\ell = \frac{\partial}{\partial v}$ covariantly constant

\subsection*{ Case 1.1.2: $\zeta = \zeta_2(u,x^e) \ell$ }

\noindent With $R_{1212} = \sigma < 0$, the metric functions $H$ and $W_i = m_i^e W_e$ will be of the form \eqref{c112metricfn} while $\zeta_2$ must satisfy the further constraint \eqref{c112zetaconstraint}. These $CSI$ spacetimes do not admit a covariantly constant vector.

\subsection*{ Case 1.2.1a :  $\zeta = n+ \frac{\sigma^*v^2}{2} \ell$}

\noindent The metric \eqref{Kundt} with $H$ and $W_i$ take the form \eqref{c121metricfnA}, $R_{1212}$ may be any value in $\mathbb{R}$.  There are no $CCNV$ spacetimes belonging to this subset of $CSI$ spacetimes.

\subsection*{ Case 1.2.1b :  $\zeta = n + \zeta_2(x^e) \ell $, $\zeta_2 \geq 0$ $\forall x^e$.}

For any $\ov{\zeta}{0}{2}(x^e) > 0,~ \forall x^e$, and any choice of $\ov{W}{1}{i}(x^e)$ satisfying \eqref{c11metric} with $\sigma = R_{1212} \leq 0$; the $CSI$ Kundt spacetime with $H$ and $W_i$ given in \eqref{c121metricfnB} will admit a timelike Killing vector. If $\ov{\zeta}{0}{2}(x^e) > 0,~ \forall x^e$, $\zeta = n$ will be a covariantly constant null vector. 

If this Killing vector is covariantly constant, $\ov{W}{1}{i} = 0$ and hence $\sigma = 0$, equation \eqref{azeta3n} and Lemma \eqref{lem:azeta0} imply $A_{ij} = 0$, and the remaining equations \eqref{azeta23} and \eqref{azeta2n} force $\ov{\zeta}{0}{2}$ to be constant. Thus $\zeta$ is the sum of the $CCNV$'s $\ell$ and $n$.  

\subsection*{ Case 1.2.2a : $\zeta = \ell + \zeta_2(u,x^e) n + \zeta_3(u) m^3$ }

For any $\zeta_3$, and a particular choice of $\zeta_2$ such that it satisfies the inequality $\zeta_3^2 \leq 2\zeta_2$ and the differential equation \eqref{c122zeta2A}, the vector $\zeta$ will be a Killing vector for the $CSI$ spacetime with metric functions given in \eqref{c122metricfnA} where $\ov{W}{1}{n}$ and $\ov{W}{0}{n}$ are, respectively, solutions \eqref{c122w1nA} and \eqref{c122w0nA}. Due to the vanishing of $\sigma^* = 4\sigma + \ov{W}{1}{i} W^{(1)i}$, the $\ov{W}{1}{i}$ must also satisfy \eqref{c11metric}

Requiring  $\zeta$ to be a $CCNV$, the $\ov{W}{1}{i}$ must vanish, causing $H = 0$ and $\sigma = 0$ as well;  this is an example of a $CCNV$ metric with $\ell = \frac{\partial}{\partial{v}}$ as the $CCNV$ and $\zeta$ acting as a second $CCNV$. The additional constraints \eqref{azeta12} - \eqref{azeta3n} require that $A_{ij} = 0$ along with the following simplification of $\zeta$ and the metric functions:

\beq & \zeta = n + \zeta_2 \ell + m^3,~~~\zeta_2 \in \mathbb{R} & \nonumber\\
& H = 0,~~~~ W_3(x^e) =  w_3(x^e),~~~~  W_n(x^e) = \int D_n (w_3) dx^3 + w_n( x^r),& \nonumber \eeq 

\noindent If $\zeta$ is timelike, then $\zeta_2 > \frac12$, while if $\zeta$ is null, $\zeta_2 = \frac12$. All of the $CSI$ spacetimes belonging to this subcase are automatically $CCNV$ with $\ell$ as another covariantly constant vector.


\subsection*{Case 1.2.2b : $n + [\frac{\sigma^*v^2}{2} + \zeta_2(u,x^e)]\ell + m^3$}

\noindent For a particular choice of $\zeta_2$ satisfying \eqref{c122zeta2B} the vector $\zeta$ is a Killing vector for the $CSI$ spacetime with the metric functions in \eqref{c122metricfnB} where the $\ov{W}{1}{n}$ and $\ov{W}{0}{n}$ satisfy \eqref{c122w1nB} and \eqref{c122w0nB}. The magnitude condition requires $\sigma^* > 0$ implying that $R_{1212} = \sigma > 0$. 
 If $\zeta$ is now covariantly constant, a contradiction arises from \eqref{azeta12}, as it requires $D_1 H = \sigma v = 0$ despite the fact that we have assumed $\sigma \neq 0$. Thus the subset of $CSI$ spacetimes associated with this subcase are never $CCNV$. 

\subsection*{Case 2 }

 Using a coordinate transform of type \eqref{type3ct} with $g(u) = \frac{u}{\sqrt{|\sigma|}}$,  $\sigma$ in equation  \eqref{HWiSigma} is rescaled  so that it equals $\sigma = -1,0,1$.  Another coordinate transform of type \eqref{type2ct} with $h = - \frac{\ov{\zeta}{0}{3}}{D_3\zeta_1}$ causes all but one component to vanish:   

\beq &W_3(u,x^3) = D_3ln(D_3\zeta_1) v,~~ W_n = 0&.  \nonumber \eeq

The other Killing vector components may be expressed entirely in terms of $\zeta_1$

\beq & \zeta_2 = \overbrace{ ( \frac{ \sigma^* \zeta_1}{4} - D_3D_3\zeta_1 ) }^{\ov{\zeta}{2}{2} } \frac{v^2}{2} + \overbrace{ ( \ov{H}{1}{~} \zeta_1  - D_2 \zeta_1 ) }^{\ov{\zeta}{1}{2} } v + \ov{\zeta}{0}{2}(u,x^e), & \nonumber \\
& \zeta_3 = - D_3 (\zeta_1) v. \nonumber \eeq

Making one final coordinate transform of type \eqref{type3ct} with $g = \int e^{-B(u)} du$ removes the $u$ dependence from $\zeta_1$ and, in fact, removes all $u$ dependence from the other components of the Killing vector and the Killing equations, (i.e., \eqref{c2-1ke22} and \eqref{c2-0ke23})  involving $\ov{H}{0}{~}$. Solving these yields the following algebraic equation for $\ov{H}{0}{~}$

\beq  \ov{H}{0}{~} = \frac{D_3( D_3(\zeta_1) \ov{\zeta}{0}{2} ) + D_3D_3ln((D_3\zeta_1)^{-\frac12}) \ov{\zeta}{0}{2} \zeta_1 }{ \zeta_1^2 D_3D_3 ln( \zeta_1 (D_3 \zeta_1)^{-\frac12})} \nonumber \eeq

\noindent With all of the Killing equations satisfied, the non-spacelike conditions \eqref{c2magv0} and \eqref{c2magv2} give two subcases depending on whether $\zeta$ is a null or timelike Killing vector. 

\subsection*{ Case 2.1: $ \zeta = \zeta_1 n + \left[ \left( \frac{\sigma^* \zeta_1}{4} - D_3D_3(\zeta_1) \right) \frac{{v}^2}{2} \right] \ell + [-D_3(\zeta_1) v] m^3 $ }

\noindent If $\zeta$ is null, $\zeta_1$ takes the following form, depending on the sign of $\sigma$:

\beq & \sigma = -1 : \zeta_1 = \frac{ (c_1 cosh(x^3) + c_2 sinh(x^3))^{-1} }{  c_1 sinh(x^3) + c_2 cosh(x^3) \pm \sqrt{ c_1^2 + c_2^2 + \left(c_1 sinh(x^3) + c_2cosh(x^3) \right)^2 } } & \nonumber \\
\nonumber \\
& \sigma = 0 : \zeta_1 = \frac{ 1}{ c_1 (1 \pm \sqrt{2}) (c_1 x^3 + c_2)} & \nonumber \\
\nonumber \\
& \sigma = 1 : \zeta_1 = \frac{(c_1 cos(x^3) + c_2 sin(x^3))^{-1}}{  -c_1 sin(x^3) + c_2 cos(x^3) \pm \sqrt{ c_1^2 + c_2^2 + \left(- c_1 sin(x^3) + c_2 cos(x^3) \right)^2 } } & \nonumber \eeq

\noindent There are no covariantly constant null vectors in this subcase as the constraint in Lemma \eqref{lem:azeta0} on $\ov{W}{1}{3}$ along with the identity \eqref{zeta3iskey} lead to a contradition with the given form of $\zeta_1$. 

\subsection*{ Case 2.2: $ \zeta = \zeta_1 n + \left[ \left( \frac{\sigma^* \zeta_1}{4} - D_3D_3(\zeta_1) \right) \frac{{v}^2}{2} + \ov{\zeta}{0}{2}(x^3) \right] \ell + [-D_3(\zeta_1) v] m^3 $}

If $\zeta$ is to be timelike, depending on the sign of $\sigma$, $\zeta_1$ takes the form: 

\beq & \sigma = -1 : \zeta_1 = \frac{sinh(x^3)}{c_1(c_1cosh(x^3) + c_2sinh(x^3))} + c_3 & \nonumber \\
\nonumber \\
& \sigma = 0 : \zeta_1 = \frac{-1}{c_1(c_1 x^3 + c_2)} + c_3 \nonumber \\
\nonumber \\
& \sigma = 1 : \zeta_1 = \frac{sin(x^3)}{c_1(c_1cos(x^3) + c_2sin(x^3))} + c_3. & \nonumber \eeq

\noindent The inequality \eqref{c2magv2} restricts the choice of $c_3$  depending on the choice of $c_1$ and $c_2$: 

\beq & \sigma = -1 : [c_1^2 + c_2^2] \zeta_1^2 - 2 \left( \frac{ c_1 sinh(x^3) + c_2 cosh(x^3)}{c_1 cosh(x^3) + c_2 sinh(x^3) } \right) \zeta_1 + \frac{1}{(c_1 cosh(x^3) + c_2 sinh(x^3))^2} < 0 & \nonumber \\
\nonumber \\
& \sigma = 0 : -c_1^2 \zeta_1^2 - 2 \left( \frac{c_1}{c_1 x^3 + c_2} \right) \zeta_1 + \frac{1}{(c_1x^3 + c_2)^2} < 0 & \nonumber \\
\nonumber \\
& \sigma = 1 : -[c_1^2 + c_2^2] \zeta_1^2 - 2 \left( \frac{ -c_1 sin(x^3) + c_2 cos(x^3) }{c_1 cos(x^3) + c_2 sin(x^3)} \right)  \zeta_1 + \frac{1}{(c_1 cos(x^3) + c_2 sin(x^3))^2} < 0 & \nonumber \eeq


\noindent There are no timelike covariantly constant vectors in $CSI$ spacetimes admitting $\zeta$ as a Killing vector.

Notice in both the null and timelike case, the value of $\sigma$ restricts the domain of $x^3$. When $\sigma = 1$ the domain of $x^3$ limited to a finite interval, $x^3 \in (x^3_0, x^3_0+ \pi)$, as the value $x^3_0 = arctan(- \frac{c_1}{c_2}) $ will cause $(D_3 \zeta_1)^{-\frac12}$ to vanish. When $\sigma = 0$, $x^3 \geq -\frac{c_2}{c_1}$ to avoid singularities. In the case with $\sigma = -1$ $x^3 > x^3_0 = arctanh(-\frac{c_1}{c_2})$ when $c_1/c_2 \leq 1 $, otherwise $\zeta_1$ is regular on the whole of the real line.

\end{document}